\documentclass[reqno,12pt,a4paper]{amsart}

\usepackage{mathrsfs}
\usepackage{amssymb}
\usepackage{amsfonts}
\usepackage{latexsym}
\usepackage{amsthm}
\usepackage{graphicx}

\def\lb{\label}

\newcommand{\er}[1]{\textrm{(\ref{#1})}}

\begin{document}


\renewcommand{\theequation}{\arabic{section}.\arabic{equation}}
\theoremstyle{plain}
\newtheorem{theorem}{\bf Theorem}[section]
\newtheorem{lemma}[theorem]{\bf Lemma}
\newtheorem{corollary}[theorem]{\bf Corollary}
\newtheorem{proposition}[theorem]{\bf Proposition}
\newtheorem{definition}[theorem]{\bf Definition}

\newtheorem{remark}[theorem]{\bf Remark}

\def\a{\alpha}  \def\cA{{\mathcal A}}     \def\bA{{\bf A}}  \def\mA{{\mathscr A}}
\def\b{\beta}   \def\cB{{\mathcal B}}     \def\bB{{\bf B}}  \def\mB{{\mathscr B}}
\def\g{\gamma}  \def\cC{{\mathcal C}}     \def\bC{{\bf C}}  \def\mC{{\mathscr C}}
\def\G{\Gamma}  \def\cD{{\mathcal D}}     \def\bD{{\bf D}}  \def\mD{{\mathscr D}}
\def\d{\delta}  \def\cE{{\mathcal E}}     \def\bE{{\bf E}}  \def\mE{{\mathscr E}}
\def\D{\Delta}  \def\cF{{\mathcal F}}     \def\bF{{\bf F}}  \def\mF{{\mathscr F}}
\def\c{\chi}    \def\cG{{\mathcal G}}     \def\bG{{\bf G}}  \def\mG{{\mathscr G}}
\def\z{\zeta}   \def\cH{{\mathcal H}}     \def\bH{{\bf H}}  \def\mH{{\mathscr H}}
\def\e{\eta}    \def\cI{{\mathcal I}}     \def\bI{{\bf I}}  \def\mI{{\mathscr I}}
\def\p{\psi}    \def\cJ{{\mathcal J}}     \def\bJ{{\bf J}}  \def\mJ{{\mathscr J}}
\def\vT{\Theta} \def\cK{{\mathcal K}}     \def\bK{{\bf K}}  \def\mK{{\mathscr K}}
\def\k{\kappa}  \def\cL{{\mathcal L}}     \def\bL{{\bf L}}  \def\mL{{\mathscr L}}
\def\l{\lambda} \def\cM{{\mathcal M}}     \def\bM{{\bf M}}  \def\mM{{\mathscr M}}
\def\L{\Lambda} \def\cN{{\mathcal N}}     \def\bN{{\bf N}}  \def\mN{{\mathscr N}}
\def\m{\mu}     \def\cO{{\mathcal O}}     \def\bO{{\bf O}}  \def\mO{{\mathscr O}}
\def\n{\nu}     \def\cP{{\mathcal P}}     \def\bP{{\bf P}}  \def\mP{{\mathscr P}}
\def\r{\rho}    \def\cQ{{\mathcal Q}}     \def\bQ{{\bf Q}}  \def\mQ{{\mathscr Q}}
\def\s{\sigma}  \def\cR{{\mathcal R}}     \def\bR{{\bf R}}  \def\mR{{\mathscr R}}
\def\S{\Sigma}  \def\cS{{\mathcal S}}     \def\bS{{\bf S}}  \def\mS{{\mathscr S}}
\def\t{\tau}    \def\cT{{\mathcal T}}     \def\bT{{\bf T}}  \def\mT{{\mathscr T}}
\def\f{\phi}    \def\cU{{\mathcal U}}     \def\bU{{\bf U}}  \def\mU{{\mathscr U}}
\def\F{\Phi}    \def\cV{{\mathcal V}}     \def\bV{{\bf V}}  \def\mV{{\mathscr V}}
\def\P{\Psi}    \def\cW{{\mathcal W}}     \def\bW{{\bf W}}  \def\mW{{\mathscr W}}
\def\o{\omega}  \def\cX{{\mathcal X}}     \def\bX{{\bf X}}  \def\mX{{\mathscr X}}
\def\x{\xi}     \def\cY{{\mathcal Y}}     \def\bY{{\bf Y}}  \def\mY{{\mathscr Y}}
\def\X{\Xi}     \def\cZ{{\mathcal Z}}     \def\bZ{{\bf Z}}  \def\mZ{{\mathscr Z}}
\def\vs{\varsigma}

\def\be{{\bf e}}
\def\bv{{\bf v}} \def\bu{{\bf u}}
\def\O{\Omega}
\def\bbD{\pmb \Delta}
\def\mm{\mathrm m}
\def\mn{\mathrm n}

\newcommand{\mc}{\mathscr {c}}

\newcommand{\gA}{\mathfrak{A}}          \newcommand{\ga}{\mathfrak{a}}
\newcommand{\gB}{\mathfrak{B}}          \newcommand{\gb}{\mathfrak{b}}
\newcommand{\gC}{\mathfrak{C}}          \newcommand{\gc}{\mathfrak{c}}
\newcommand{\gD}{\mathfrak{D}}          \newcommand{\gd}{\mathfrak{d}}
\newcommand{\gE}{\mathfrak{E}}
\newcommand{\gF}{\mathfrak{F}}           \newcommand{\gf}{\mathfrak{f}}
\newcommand{\gG}{\mathfrak{G}}           
\newcommand{\gH}{\mathfrak{H}}           \newcommand{\gh}{\mathfrak{h}}
\newcommand{\gI}{\mathfrak{I}}           \newcommand{\gi}{\mathfrak{i}}
\newcommand{\gJ}{\mathfrak{J}}           \newcommand{\gj}{\mathfrak{j}}
\newcommand{\gK}{\mathfrak{K}}            \newcommand{\gk}{\mathfrak{k}}
\newcommand{\gL}{\mathfrak{L}}            \newcommand{\gl}{\mathfrak{l}}
\newcommand{\gM}{\mathfrak{M}}            \newcommand{\gm}{\mathfrak{m}}
\newcommand{\gN}{\mathfrak{N}}            \newcommand{\gn}{\mathfrak{n}}
\newcommand{\gO}{\mathfrak{O}}
\newcommand{\gP}{\mathfrak{P}}             \newcommand{\gp}{\mathfrak{p}}
\newcommand{\gQ}{\mathfrak{Q}}             \newcommand{\gq}{\mathfrak{q}}
\newcommand{\gR}{\mathfrak{R}}             \newcommand{\gr}{\mathfrak{r}}
\newcommand{\gS}{\mathfrak{S}}              \newcommand{\gs}{\mathfrak{s}}
\newcommand{\gT}{\mathfrak{T}}             \newcommand{\gt}{\mathfrak{t}}
\newcommand{\gU}{\mathfrak{U}}             \newcommand{\gu}{\mathfrak{u}}
\newcommand{\gV}{\mathfrak{V}}             \newcommand{\gv}{\mathfrak{v}}
\newcommand{\gW}{\mathfrak{W}}             \newcommand{\gw}{\mathfrak{w}}
\newcommand{\gX}{\mathfrak{X}}               \newcommand{\gx}{\mathfrak{x}}
\newcommand{\gY}{\mathfrak{Y}}              \newcommand{\gy}{\mathfrak{y}}
\newcommand{\gZ}{\mathfrak{Z}}             \newcommand{\gz}{\mathfrak{z}}

\def\ve{\varepsilon}   \def\vt{\vartheta}    \def\vp{\varphi}    \def\vk{\varkappa}

\def\A{{\mathbb A}} \def\B{{\mathbb B}} \def\C{{\mathbb C}}
\def\dD{{\mathbb D}} \def\E{{\mathbb E}} \def\dF{{\mathbb F}} \def\dG{{\mathbb G}}
\def\H{{\mathbb H}}\def\I{{\mathbb I}} \def\J{{\mathbb J}} \def\K{{\mathbb K}}
\def\dL{{\mathbb L}}\def\M{{\mathbb M}} \def\N{{\mathbb N}} \def\dO{{\mathbb O}}
\def\dP{{\mathbb P}} \def\R{{\mathbb R}}\def\S{{\mathbb S}} \def\T{{\mathbb T}}
\def\U{{\mathbb U}} \def\V{{\mathbb V}}\def\W{{\mathbb W}} \def\X{{\mathbb X}}
\def\Y{{\mathbb Y}} \def\Z{{\mathbb Z}}


\def\la{\leftarrow}              \def\ra{\rightarrow}            \def\Ra{\Rightarrow}
\def\ua{\uparrow}                \def\da{\downarrow}
\def\lra{\leftrightarrow}        \def\Lra{\Leftrightarrow}


\def\lt{\biggl}                  \def\rt{\biggr}
\def\ol{\overline}               \def\wt{\widetilde}
\def\ul{\underline}
\def\no{\noindent}


\let\ge\geqslant                 \let\le\leqslant
\def\lan{\langle}                \def\ran{\rangle}
\def\/{\over}                    \def\iy{\infty}
\def\sm{\setminus}               \def\es{\emptyset}
\def\ss{\subset}                 \def\ts{\times}
\def\pa{\partial}                \def\os{\oplus}
\def\om{\ominus}                 \def\ev{\equiv}
\def\iint{\int\!\!\!\int}        \def\iintt{\mathop{\int\!\!\int\!\!\dots\!\!\int}\limits}
\def\el2{\ell^{\,2}}             \def\1{1\!\!1}
\def\sh{\sharp}
\def\wh{\widehat}
\def\bs{\backslash}
\def\intl{\int\limits}

\def\na{\mathop{\mathrm{\nabla}}\nolimits}
\def\sh{\mathop{\mathrm{sh}}\nolimits}
\def\ch{\mathop{\mathrm{CR}}\nolimits}
\def\where{\mathop{\mathrm{where}}\nolimits}
\def\all{\mathop{\mathrm{all}}\nolimits}
\def\as{\mathop{\mathrm{as}}\nolimits}
\def\Area{\mathop{\mathrm{Area}}\nolimits}
\def\arg{\mathop{\mathrm{arc}}\nolimits}
\def\const{\mathop{\mathrm{const}}\nolimits}
\def\det{\mathop{\mathrm{det}}\nolimits}
\def\diag{\mathop{\mathrm{diag}}\nolimits}
\def\diam{\mathop{\mathrm{diam}}\nolimits}
\def\dim{\mathop{\mathrm{dim}}\nolimits}
\def\dist{\mathop{\mathrm{dist}}\nolimits}
\def\Im{\mathop{\mathrm{Im}}\nolimits}
\def\Iso{\mathop{\mathrm{Iso}}\nolimits}
\def\Ker{\mathop{\mathrm{Ker}}\nolimits}
\def\Lip{\mathop{\mathrm{Lip}}\nolimits}
\def\rank{\mathop{\mathrm{rank}}\limits}
\def\Ran{\mathop{\mathrm{Ran}}\nolimits}
\def\Re{\mathop{\mathrm{Re}}\nolimits}
\def\Res{\mathop{\mathrm{Res}}\nolimits}
\def\res{\mathop{\mathrm{res}}\limits}
\def\sign{\mathop{\mathrm{sign}}\nolimits}
\def\span{\mathop{\mathrm{span}}\nolimits}
\def\supp{\mathop{\mathrm{supp}}\nolimits}
\def\Tr{\mathop{\mathrm{Tr}}\nolimits}
\def\BBox{\hspace{1mm}\vrule height6pt width5.5pt depth0pt \hspace{6pt}}


\newcommand\nh[2]{\widehat{#1}\vphantom{#1}^{(#2)}}
\def\dia{\diamond}

\def\Oplus{\bigoplus\nolimits}



\def\qqq{\qquad}
\def\qq{\quad}
\let\ge\geqslant
\let\le\leqslant
\let\geq\geqslant
\let\leq\leqslant
\newcommand{\ca}{\begin{cases}}
\newcommand{\ac}{\end{cases}}
\newcommand{\ma}{\begin{pmatrix}}
\newcommand{\am}{\end{pmatrix}}
\renewcommand{\[}{\begin{equation}}
\renewcommand{\]}{\end{equation}}
\def\eq{\begin{equation}}
\def\qe{\end{equation}}
\def\[{\begin{equation}}
\def\bu{\bullet}

\title[{Eigenvalues of periodic difference operators on lattice octant}]
{Eigenvalues of Periodic difference operators on lattice octant}

\date{\today}
\author[Evgeny Korotyaev]{Evgeny Korotyaev}
\address{Saint-Petersburg State University, Universitetskaya nab. 7/9, St. Petersburg, 199034, Russia,
\ korotyaev@gmail.com, \ e.korotyaev@spbu.ru,}


\begin{abstract}

Consider a difference operator $H$ with periodic coefficients on the
 octant of the lattice. We show that for any integer $N$ and any bounded
interval $I$, there exists an operator $H$ having $N$ eigenvalues,
counted  with multiplicity on this interval, and does not exist
†other spectra on the interval. Also right and to the left of it are
spectra and the corresponding subspaces have an infinite dimension.
 Moreover, we prove similar results for other domains and any dimension. The
proof is based on the inverse spectral theory for periodic Jacobi
operators.

\end{abstract}

\maketitle

\vskip 0.25cm

\section {\lb{Sec0}Introduction and main results}

We consider a operator  $H=H_1+H_2+V$ acting on domain $D$ and $V$
is the multiplication operator on $\ell^2(D)$:
\[
\lb{d1}
\begin{aligned}
(Vf)_z=V(z)f_z, \qqq f=f_z, \qq z=(x,y)\in D=\Z_+^{d_1}\ts
\Z^{d_2}\ss \Z^d,
\end{aligned}
\]
$d_1+d_2=d\ge 2,\  d_1,d_2\ge  0$.
 Here  $H_1$ is the difference
operator on the octant $ \Z_+^{d_1}$ with the Dirichlet boundary
conditions on the boundary $\pa \Z_+^{d_1}$ (i.e., $g=0$ on
$\Z^{d_1}\sm \Z_+^{d_1}$ in \er{d2}) and $H_2$ is the difference
operator on $\Z^{d_2}$ defined by
\[
\lb{d2}
\begin{aligned}
(H_1 g)_x =\sum_{i=1}^{d_1}
(a_{x-\vs_i}^ig_{x-\vs_i}+a_x^ig_{x+\vs_i}), \qq x\in \Z_+^{d_1}, \qqq g=(g_x)_{}\in \ell^2(\Z_+^{d_1}),\\
(H_2 u)_y =\sum_{i=d_1+1}^{d} (a_{y-\vs_i}^i
u_{y-\vs_i}+a_y^iu_{y+\vs_i}),\qq y\in \Z^{d_2}, \qqq u=(u_y)\in
\ell^2(\Z^{d_2}).
\end{aligned}
\]
Here $ \vs_1=(1,0,0,..),..., \vs_d=(0,0,0,..,1)$ is the standard
basis in $\Z^d$ and $\Z_+=\{1,2,3,...\}$. We assume that the
potential $V$ and the coefficients $a^i$ are real octant periodic,
i.e., they have decompositions \er{VO}. In order to define octant
periodic functions we introduce a sequence $\o=(\o_j)_1^m$, where
$\o_j=+$ or $\o_j=-$ and the set of all such sequences we denote by
$\O_m$. For any $\o\in \O_m$ we define the octants $\cZ_\o\ss
\Z^{m}$ by
$$
\cZ_\o=\Z_{\o_1}\ts \Z_{\o_2}\ts ....\ts
\Z_{\o_m},\qqq\o=(\o_j)_1^m\in \O_m, \qq \Z_-=\Z\sm
\Z_+=\{...,-3,,-2,-1,0\}.
$$
In particular, if  $d=2$ and $\o=(+,+)$, then we have the positive
octant $\cZ_\o=\Z_+^2$. Note that two axial lines $(x_1={1\/2},
x_2={1\/2})$ divide space $\Z^2$ into four quadrants, each with a
coordinate signs from $(-,-)$ to $(+,+)$.

{\it  A function $F(z), z=(x,y)\in \Z_+^{d_1}\ts \Z^{d_2}$ is called
octant periodic if it   has the decomposition
\[
\lb{VO}
\begin{aligned}
  F(x,y)=\sum_{\o\in \O_{d_2}} F_\o(x,y)\c_\o(y),\qqq
\end{aligned}
\]
where $\c_\o$ is the characteristic function of the octant $\cZ_\o$
and the function $F_\o(z), z=(x,y)$ is periodic in $\Z^d$ and
satisfies
\[
\lb{V}
\begin{aligned}
 F_\o(z+\gp_i \vs_i)=F_\o(z),\qq {\rm for \  all} \ (z,i)\in \Z^d\ts \N_d
\end{aligned}
\]
for some constants $\gp_i=\gp_i(\o)>0$, where $\N_d=\{1,2,...,d\}$.
}

For each $\o\in\O_d$ we define difference operators $H_\o$ with
periodic functions $a^{\o,i},  V_\o$ on $\Z^d$  by
\[
\lb{dHo}
\begin{aligned}
 H_\o f=\sum_{i=1}^{d} ((a_z^{\o,i}f_{z+\vs_i}
+a_{z-\vs_i}^{\o,i}f_{z-\vs_i})) +V_\o, \qqq f=(f_z)\in
\ell^2(\Z^{d}).
\end{aligned}
\]

It is well known that the spectrum of each operator $H_\o$ is
absolutely continuous and  is an union of  a finite number of
bounded intervals. In the next theorem we show the existence of
eigenvalues of $H$ with some octant periodic functions $a_z^{\o,i} $
and $ V_\o$.

\begin{theorem}\lb{T1}
i) Let an operator $H$ be given by \er{d1}, \er{d2} with octant
periodic coefficients. Then
\[
\lb{spH} \bigcup_{\o\in \O_{d_2}}\s(H_\o) \subseteq \s_{ess}(H).
\]
ii)  Let $I\ss \R$ be a finite open interval. Then for any integer
$N\ge 0$ there exists an operator $H$ given by \er{d1}, \er{d2} and
having $N$ eigenvalues, counted  with multiplicity on this interval,
and does not exist †other spectra on the interval. Also right and to
the left of it are spectra and the corresponding subspaces have an
infinite dimension.
\end{theorem}

\no \textbf{Remark.} 1) The result of i) is standard and its proof
is based on the Floquet theory.

2) We do not know any information about absolutely continuous
spectrum of $H$. We only show  that the operator $H$ can have any
number $N\ge 1$ of eigenvalues for specific coefficients.

3) In the  case of the continuous Schr\"odinger operator with octant
periodic potentials on $\R^d$
 the top of the spectra is a isolated simple eigenvalue for
specific potentials \cite{KM18}. In the discrete case we have no any
information about it.

\subsection{Historical review}

Firstly we discuss the continuous case. The one dimensional model of
octant periodic potentials on $\R^1$ is considered by Korotyaev
\cite{K00}, \cite{K05}.  The corresponding multidimensional model of
octant periodic potentials is considered recently by Korotyaev and
Moller \cite{KM18}. Hempel and Kohlmann \cite{HK11},\cite{HK11x}
discuss different types of dislocation problem in solid state
physics.

Secondly, we discuss the discrete case. Local defects are considered
by different authors. For the discrete Schr\"odinger operators most
of the results were obtained for uniformly decaying potentials for
the $\Z$ case, see, for example, \cite{T89}. There are results about
spectral properties of discrete Schr\"odinger  operators on the
lattice $\Z^d$, the simplest example of periodic graphs.
Schr\"odinger operators with decaying potentials on the lattice
$\Z^d$ are considered by Boutet de Monvel-Sahbani \cite{BS99},
Hundertmark-Simon \cite{HS02}, Isozaki-Korotyaev \cite{IsK12},
Isozaki-Morioka \cite{IM14}, Korotyaev-Moller \cite{KM17}, Nakamura
\cite{Na14}, Parra and Richard \cite{PR18}, Rosenblum-Solomjak
\cite{RoS09}, Shaban-Vainberg \cite{SV01}  and see references
therein. Gieseker-Kn\"orrer-Trubowitz \cite{GKT93} consider
Schr\"odinger operators with periodic potentials on the lattice
$\Z^2$. Korotyaev-Kutsenko \cite{KK10} study the spectra of the
discrete Schr\"odinger operators on graphene  nano-ribbons in
external electric  fields. The inverse spectral theory for the
discrete Schr\"odinger operators with finitely supported potentials
on some graphs were discussed by
  Ando \cite{A13}, Ando-Isozaki-Morioka
\cite{AIM16}, \cite{AIM18}, Isozaki-Korotyaev \cite{IsK12}.
Scattering on periodic metric graphs was considered by
Korotyaev-Saburova \cite{KS15}.
 Laplacians on periodic
graphs with non-compact perturbations and the stability of their
essential spectrum were considered in \cite{HKSV15}, \cite{SS15}.
Korotyaev-Saburova \cite{KS17},  \cite{KS17x} considered
Schr\"odinger operators with periodic potentials on periodic
discrete graphs with by so-called guides, which are periodic in some
directions and finitely supported in others. They described some
properties of so-called guided spectrum.  Note that line defects on
the lattice were considered in \cite{C12}, \cite{Ku14}, \cite{Ku16},
\cite{OA12}.     Hempel, Kohlmann, Stautz and Voigt \cite{HKSV15}
discussed  nano-tubes     with a dislocation.

 We shortly describe the plan of the paper. In Section~2 we present the
main properties of the periodic Jacobi operator on the half lattice
$\Z_+$. In Section~3 we discuss half-solid models on the lattice
$\Z$. Section~4 is a collection of needed facts about difference
operators on $D$, when the variables are separated. In Section~5 we
prove main theorems.


\section{Periodic Jacobi operators on the half-lattice}
\setcounter{equation}{0}

\subsection{Periodic Jacobi operators}
Let $\N_\pm=\{\pm 1, \pm2,\pm3,...\}$. Recall that
$\N_+=\Z_+=\{1,2,3,...\}$ and $\Z_-=\Z\sm \Z_+$.  We consider the
p-periodic Jacobi operator
 $J_{_\pm }$ on $\ell^2(\N_\pm )$ given by
\[
\begin{aligned}
\label{Jdef} (J_\pm f)_x=a_{x-1}f_{x-1}+a_{x}f_{x+1}+b_xf_x,\qqq
f_0=0,\qq x\in \N_\pm,
\end{aligned}
\]
and in particular,
$$
\begin{aligned}
\ca (J_{_-} yf)_{-1}=a_{-3}f_{-3}+a_{-2}f_{1}+b_{-2}f_{-2}, \\
 (J_{_-} yf)_{-1}=a_{-2}f_{-2}+b_{-1}f_{-1},\qq f_1=0 \ac,
\qqq
\ca (J_{_+} f)_2=a_{1}f_{1}+a_{2}f_{3}+b_2f_2, \\
 (J_{_+}f)_1=a_{1}f_{2}+b_1f_1\qq f_0=0  \ac,
\end{aligned}
$$
where $a_x>0, b_x\in\R, x\in \Z$ are $p$ periodic sequences  and the
product $\prod_1^p a_j=1$.
 It is well known that the spectrum of $J_{_\pm}$ has absolutely
continuous part $\s_{ac}(J_{_-})=\s_{ac}(J_{_+})$ (the union of the
bands $\s_0, \s_n, n\in\N_{p-1}$ separated by gaps $\g_n$) plus at
most one eigenvalue of $J_{_-}$ or $J_{_+}$ in each non-empty gap
$\g_x$, $n\in\N_{p-1}$. The bands $\s_n$ and gaps $\g_n$ are given
by
\[
\begin{aligned}
\label{DA2}
& \s_0=[\l_{0}^+,\l_{1}^-], \qq
\s_n=[\l_{n}^+,\l_{n+1}^-], \qqq \g_{n}=(\l^-_{n},\l^+_n),\qq n\in
\N_{p-1},
\\
&   \l_0^{+}<\l_1^-\le \l_1^+< ... < \l_{p-1}^- \le \l_{p-1}^+ <
\l_{p}^{-},
 \end{aligned}
 \]
 (see
Fig. 1)
and recall that $\N_p=\{1,2,...,p\}$. The bands satisfy (see e.g.,
\cite{L92}, \cite{KKr03})
\[
\lb{eb} \sum_{n=0}^{p-1} |\s_n|\le 4.
\]
\begin{figure}
\tiny \unitlength=1.00mm \special{em:linewidth 0.4pt}
\linethickness{0.4pt}
\begin{picture}(108.67,33.67)
\put(41.00,17.33){\line(1,0){67.67}}
\put(44.33,9.00){\line(0,1){24.67}}
\put(108.33,14.00){\makebox(0,0)[cc]{$\Re\l$}}
\put(41.66,33.67){\makebox(0,0)[cc]{$\Im\l$}}
\put(44.33,17.33){\linethickness{4.0pt}\line(1,0){11.33}}
\put(66.66,17.33){\linethickness{4.0pt}\line(1,0){11.67}}
\put(82.00,17.33){\linethickness{4.0pt}\line(1,0){12.00}}
\put(95.66,17.33){\linethickness{4.0pt}\line(1,0){11.00}}
\put(46.66,20.00){\makebox(0,0)[cc]{$\l_0^+$}}
\put(56.66,20.33){\makebox(0,0)[cc]{$\l_1^-$}}
\put(68.66,20.33){\makebox(0,0)[cc]{$\l_1^+$}}
\put(78.33,20.33){\makebox(0,0)[cc]{$\l_2^-$}}
\put(84.33,20.33){\makebox(0,0)[cc]{$\l_2^+$}}
\put(93.00,20.33){\makebox(0,0)[cc]{$\l_3^-$}}
\put(98.66,20.33){\makebox(0,0)[cc]{$\l_3^+$}}
\put(106.33,20.33){\makebox(0,0)[cc]{$\l_4^-$}}
\end{picture}
\caption{  \footnotesize The cut domain $\C\sm \cup \s_n$ and the
cuts (bands) $\s_n=[\l^+_{n},\l^-_{n+1}], n=0,1,...$} \lb{sS}
\end{figure}
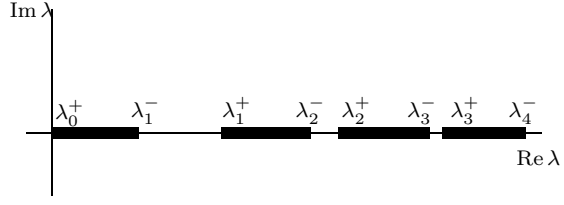
If a gap $\g_n$ is degenerate, i.e. $|\g_n|=0$, then the
corresponding segments $\s_n$, $\s_{n+1}$ merge. We  introduce
fundamental solutions $\vp=(\vp_x(\l))_{x\in\Z}$ and $
\vt=(\vt_x(\l))_{x\in \Z} $ of the equation
\[
 \lb{b}
a_{x-1} f_{x-1}+a_{x}f_{x+1}+b_xf_x=\l f_x,\ \ (\l,x)\in\C\ts\Z,
\]
with initial conditions $\vp_{0}=\vt_1=0$ and $\vp_1=\vt_{0}=1$.
Recall that the zeros of $\vp_x(\l)$ are real, simple and strictly
interlace those of $\vp_{x+1}(\l)$. Moreover, the zeros of
$\vt_x(\l)$ are real, simple and strictly interlace those of
$\vp_{x}(\l)$. Define the Lyapunov function $\gF$ by
$$
\gF(\l)={1\/2}(\vp_{p+1}(\l)+\vt_p(\l)), \qqq \l\in \C.
$$
We recall the well known asymptotics as $\l\to \iy$:
\[
\lb{1} \vt_x(\l)=-{a_0\l^{x-2} (1+O({1\/\l}))\/a_1..a_{x-1}},\qqq
\vp_x(\l)={\l^{x-1}(1+O({1\/\l}))\/a_1..a_{x-1}},
\]
\[
\lb{2} \gF(\l)={\l^{p}\/2}+O(\l^{p-2})....
\]
The functions $\gF, \vp_x$ and $\vt_x, x\ge1$ are polynomials of
$(\l,a,b)\in\C^{2p+1}$. We have the following identities
$$
\s_{ac}(J_{\pm})=\{\l\in\R:\ |\gF(\l)|\le 1\}\qqq  and \qqq
(-1)^{p-n}\gF(\l_n^{\pm})=1,\qq n=0,1,2,....,p.
$$
 For any sequences
$u=(u_x)_{-\iy}^{\iy}, f=(f_x)_{-\iy}^{\iy}$ we define the Wronskian
\[
\{f,u\}_x=a_x(f_xu_{x+1}-u_xf_{x+1}),\qq x\in \Z.
\]
If $f,u$ are some solutions of \er{b}, then $\{f,u\}_x$ does not
depend on $x$. In particular, we have
\[
\lb{wfs} \vt_p \vp_{p+1}-\vp_p\vt_{p+1}=1,
\]
since $\{\vt,\vp\}_p=a_p(\vt_p
\vp_{p+1}-\vp_p\vt_{p+1})=\{\vt,\vp\}_0=a_0$. Thus we obtain
\[
\lb{wf} \gF_o^2-\gF^2+1=-\vp_p\vt_{p+1}.
\]

We define the Jacobi operator $J_D$ on $\N_{p-1}$  with the
Dirichlet boundary conditions by
\[
\lb{D} (J_D f)_x=a_{x-1}f_{x-1}+a_{x}f_{x+1}+b_xf_x,\qqq x\in
\N_{p-1}, \qq f_0=f_p=0.
\]
  Denote its corresponding Dirichlet eigenvalues by  $\m_n, n\in \N_{p-1}$. It is
well known that the eigenvalues $\m_n$ are simple and  are zeros of
the polynomial $\vp_p(\l)$ and satisfy
\[
\lb{Np1}  \m_n\in [\l_n^-, \l_n^+]\qq {\rm for \ all} \ n\in
\N_{p-1}.
\]

\subsection{Riemann surface}
For the  operator $J_\pm$ we  introduce the two-sheeted Riemann
surface $\L$ obtained by joining the upper and lower rims of two
copies of the cut plane $\C\sm\s_{ac}(J_+)$ in the usual (crosswise)
way. We denote the $n$-th gap on the first physical sheet $\L_1$ by
$\g_n^{1}$ and the same gap but on the second nonphysical sheet
$\L_2$ by $\g_n^{2}$, and set a circle gap $\g_n^c$ by
$$
\g_n^c =\ol\g_n^{1}\cup \ol\g_n^{2}\qq {\rm for \ all} \ n\in
\N_{p-1},
$$
see Fig. \ref{fig2}.
 Note that $\L$ is the two-sheeted
Riemann surface for $\sqrt{1-\gF^2(\l)}$.  The polynomial  $\gF(\l)$
is real on the real line.
 We use the standard definition of the root: $\sqrt{1}=1$ and fix the branch of the function
 $\gf(\l)=\sqrt{1-\gF^2(\l)}$ on $\C$ by demanding
\[
\lb{rs}
 \gf(\l)=\sqrt{1-\gF^2(\l)} <0,\qqq \mbox{for} \qq \l\in (\l_{p-1}^+, \l_{p}^-).
\]

\setlength{\unitlength}{1.0mm}
\begin{figure}[h]
\centering
\unitlength 1.0mm 
\begin{picture}(90,30)

\bezier{600}(30,0)(55.5,3)(81,6) \bezier{600}(34,6)(59.5,9)(85,12)

\bezier{600}(18,15)(28,14)(34,6) \bezier{600}(69,21)(79,20)(85,12)

\bezier{600}(18,15)(20,6)(30,0) \bezier{35}(69,21)(71,12)(81,6)
\bezier{100}(74.8,10.7)(76.2,9.2)(81,6)

\bezier{600}(18,15)(28.5,16.2)(39,17.4)
\bezier{600}(48,18.5)(58.5,19.75)(69,21) \put(39,17.4){\circle*{1}}
\put(48,18.5){\circle*{1}} \bezier{20}(39,17.4)(44.5,13)(48,18.5)
\bezier{20}(39,17.4)(43.0,23.3)(48,18.5)

\put(37,14){$\scriptstyle\l_1^-$}
\put(48,15.5){$\scriptstyle\l_1^+$}

\bezier{600}(18,15)(16,20)(4,24) \bezier{600}(69,21)(67,26)(55,30)

\bezier{600}(18,15)(7,15)(0,18) \bezier{35}(69,21)(58,21)(51,24)

\bezier{100}(0,18)(25.5,21)(51,24)
\bezier{100}(0,18)(6.9,18.85)(13.8,19.7)

\bezier{600}(4,24)(29.5,27)(55,30)

\end{picture}

\caption{\footnotesize The two-sheeted Riemann surface with the open
gap $(\l_1^-,\l_1^+)$ and the circle gap} \label{fig2}
\end{figure}
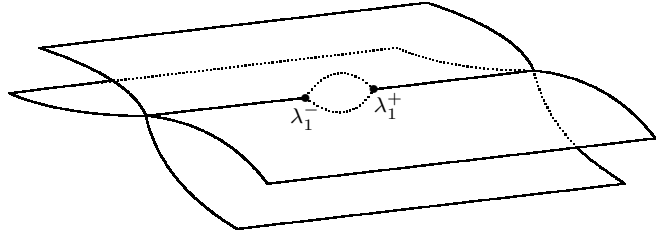

\subsection{Bloch functions}

 Define the cut spectral domain $\dL$ and the cut quasimomentum domain
$\K_p$ by
\[
\label{Km1}
\begin{aligned}
& \dL=\C\sm \cup_{n=1}^{p-1} \g_n,\qqq   \K_p=\rt\{-\pi p<\Re k<0\rt\} \sm \cup_{n=1}^{p-1} \G_n,\\
& \G_n=(\pi n+ih_n, \pi n-ih_n),\qqq {\rm ch}\ h_n=\max_{\l\in \g_n}
|\gF(\l)|=(-1)^{n-p}\gF(\a_n),
  \end{aligned}
\]
where $h_n\ge 0$ is defined by the equation $\cosh h_n=|\gF(\a_n)|$,
where $\a_n$ is a zero of $\gF'(\l)$ in the close gap
$[\l_n^-,\l_n^+]$. For each Jacobi operator $J_\pm$ there exist a
unique conformal mapping $k:\dL\to \K_p$ such that $\gF(\l)=\cos
k(\l),\  \l\in \dL$ and  following identities and asymptotics hold
true:
\[
\label{Km2}
\begin{aligned}
& k(\dL)=\K_p,\qqq   k(\C_\pm)=\K_p\cap \C_\pm,\qqq       k(\g_n\pm i0)=\G_n\cap \C_\pm,\qq \\
&  k(\l)\to \pm i \iy \qqq \as \ \Im \l\to \pm \iy,
  \end{aligned}
\]
see \cite{P84}.
 The quasimomentum $k(\l)$ satisfies $\sigma_{ac}(J_\pm)=\{\l\in\R;\,\,\Im k(\l)=0\}$.
 We define the Bloch functions functions $\p_x^\pm$  and the Weyl-Titchmarsh function $m_\pm$  by
\[
\label{bf}
\begin{aligned}
 \p_x^\pm (\l)=\vt_x(\l)+m_\pm(\l)\vp_x(\l),\qqq
 m_\pm(\l)={\gF_o(\l)\pm i \sin k(\l) \/ \vp_p(\l)},\qqq \l\in \L_1,
  \end{aligned}
\]
where $\sin k(\l)=\sqrt{1-\gF^2(\l)}, \l\in \L_1$ and satisfies
\er{rs}. Due to the properties of $\sin k(\l)=\sqrt{1-\gF^2(\l)},
\l\in \L_1$ the functions $\p_x^\pm$  and $m_\pm$ have  analytic
extensions from the first sheet $\L_1$ onto the whole two-sheeted
Riemann surface $\L$. Let $k=u(\l)+iv(\l)$.

\subsection{Eigenvalues and resonances}
It is well known (see e.g. \cite{IK11})  that, for each finitely
supported $f\in \ell^2(\N_\pm), f\ne 0$, the function
$g(\l)=((J_\pm-\l)^{-1}f,f)$ has a meromorphic extension from the
physical sheet $\L_1$ into the whole Riemann surface $\L$. It is
well known that the function $g(\l)$ has only tree following kinds
of singularity on $\L$:

$\bu$ $g$ has a pole at some $\l_o\in \bigcup\limits_{n=1}^{p-1}
\g_n^{(1)}\ss \L_1$ for some $f$ and $\l_o$ is an eigenvalue  of
$J_\pm$.

$\bu$ $g$ has a pole at some $\l_o\in \bigcup\limits_{n=1}^{p-1}
\g_n^{(2)}\ss \L_2$ for some $f$ and $\l_o$ is called a resonance of
$J_\pm$.

$\bu$ The function $g(\l_o +z^2)$ has a pole at $z=0$ for some
$\l_o\in\{\l_n^+, \l_n^-\}, n=\N_{p-1}$ and  $\l_o$ is  called  a
virtual state  of $J_\pm$.

We call $\l_o$ a state if $\l_o$ is an eigenvalue or a resonance or
a virtual state. It is well known that if some gap $\g_n\ne \es,
n\in \N_{p-1}$, see e.g., \cite{IK11},  then the operator $J_\pm$
has exactly one state $\m_n^\pm$ on each "circle" gap $\g_n^c$ and
there are no others. The projection of $\m_n^\pm\in \g_n^c$ onto the
complex plane coincides with the eigenvalue $\m_n$ of the operator
$J_D$ with the Dirichlet boundary conditions \er{D}.
 There are no other states of $J_\pm$. If there are exactly
$n_\bu\ge 1$ non-degenerate finite gaps in the spectrum of
$\s_{ac}(J_\pm)$, then the operator $J_\pm$ has exactly $n_\bu$
states; the closed gaps $\g_n=\es$ and the semi-infinite gaps
$(-\iy,\l_0^+)$ and $(\l_p^-,\iy,)$ do not contribute any states. In
particular, if $\g_n=\es$ for all $n\in \N_{p-1}$, then all $a_n=1,
b_n=0$ (see e.g., \cite{K06}) and thus $J_\pm$ has no states. A more
detailed description of the states of $J_\pm$ is given below.


\begin{lemma}\lb{Tbf1}
 Let $\l\in \L_1\sm \{\m_j, j\in \N_{p-1}\}$. Then
\[
\lb{bf1} m_\pm(\l)={e^{\pm
ik(\l)}-\vt_{p}(\l)\/\vp_p(\l)}={\vp_{p+1}(\l)-e^{\mp
ik(\l)}\/\vp_{p}(\l)},
\]
\[
\lb{mm} m_+m_-=-\frac{\vt_{p+1}}{\vp_p},
\]
and
\[
\lb{bf2}
\begin{aligned}
\p_0^\pm (\l)=1,\qqq \p_p^\pm (\l)=e^{\pm ik(\l)},
\\
\p_1^\pm(\l)=m_\pm(\l),\qqq \p_{p+1}^\pm(\l)=e^{\pm
ik(\l)}m_\pm(\l).
\end{aligned}
\]
\end{lemma}

{\bf Remark.} This relation \er{mm} considered at zeros of
$\vp_p(\l)$ shows that  if $\mu_n^\pm\in \g_n^c\ne \es,$ $n\in
\N_{p-1},$ is not a virtual state then we have:

\no $\bu$ $m_\pm$ has simple pole at $\mu_n^\pm$ and the function
$m_\mp$ is regular at $\mu_n^\pm$.

\no $\bu$   The solutions $\psi_n^\pm(\l)$ is regular at $\mu_n^\pm$
iff the other $\psi_n^\mp(\l)$ has simple poles at $\mu_n^\pm\in
\g_n^c$.

\no {\bf Proof}. Let for shortness $\vt_x=\vt_x(\cdot),
\vp_x=\vp_x(\cdot),...$ We consider the case $+$, the proof for the
case $-$ is similar. Using the definitions
$\gF={\vp_{p+1}+\vt_p\/2}=\cos k$ and $\gF_o={\vp_{p+1}-\vt_p\/2}$
we obtain
$$
\begin{aligned}
e^{ik}-\vt_p=\gF+i\sin k-\vt_p=\gF_o+i\sin k,
\\
\vp_{p+1}-e^{- ik}=\vp_{p+1}-\gF+i\sin k=\gF_o+i\sin k,
\end{aligned}
$$
which yields \er{bf1}. From \er{wf}, \er{wfs}  we obtain
$$
\vp_p^2 m_+m_-= \gF_o^2-(\gF^2-1)=1+\vt_p\vp_{p+1}=-\vp_p\vt_{p+1},
$$
which yields \er{mm}.
 We show \er{bf2}. From \er{b} at $n=0$ and
$n=p$  and \er{bf1}  we have
\[
\begin{aligned}
\lb{bfb} \p_0^{+}=1,\qqq \p_1^{+}=\vt_1+m_\pm\vp_1
=m_+,\qqq\p_p^{+}=\vt_p+(e^{ik}-\vt_p)=e^{ik}.
\end{aligned}
\]
 From  \er{bf1}, \er{wfs} we have
$$
\begin{aligned}
 \p_{p+1}^+=\vt_{p+1}+{e^{ ik}-\vt_p\/\vp_p }\vp_{p+1}
={\vt_{p+1}\vp_p+(e^{ik}-\vt_p)\vp_{p+1}\/\vp_p }
 =
{\vp_{p+1}e^{ik}-1\/\vp_p }=
e^{ik}m_+,
\end{aligned}
$$
which  yields \er{bf2}.
 \BBox

We recall the results from \cite{IK11}

\begin{lemma}\lb{TAn}
Let a finite gap $\g_n\ne \es$ for some $n\in\N_{p-1}$. Then

\no i) the operator $J_\pm $ has exactly one  state $\m_n^\pm$ on
$\g^{\rm c}_n= \ol\g^{1}_n\cup \ol\g^2_n$ and its  projection  on
$\C$ coincides with the Dirichlet eigenvalue $\mu_n$.

\no ii) the state  $\m_j^\pm \in \g^{1}_n$ iff the state
 $\m_j^\mp \in \g^{2}_n$. Moreover, if $\m_n^\pm \in\{\l_n^+, \l_n^-\}$ is a virtual state, then
$\m_n^+=\m_n^-$.

\no iii) Let $\m_j^+ \in \g^{1}_n$ be an eigenvalue of $J_+$. Then
$\vp(\m_j^+)\in \ell^2(\Z_+)$.

\no iv)  Let $\l_o \in \g_n^c \neq \es$ be a state of $J_+$ for some
$n \in \N_{p-1}$. Then
\[
\begin{aligned}
\lb{bfc} \ca {\rm if}\qq |\vp_{p+1}(\l_o)| < 1 \Rightarrow \ \l_n
\in \g_n^{1},
\\
{\rm if}\qq |\vp_{p+1}(\l_o)| > 1 \Rightarrow \ \l_o \in \g_n^{2},
\\
{\rm if}\qq |\vp_{p+1}(\l_o)|=1 \Rightarrow  \ \l_o \in \{\l_n^+,
\l_n^- \} \ac
\end{aligned}
\]

\end{lemma}

\no {\bf Proof}. The proof of i)-iii) is standard (see e.g.
\cite{IK11}).

iv) Let $\m_n\in \{\l_n^-, \l_n^+\}$. The Lyapunov function
$\vp_{p+1}(\m_n)+\vt_{p}(\m_n)=(-1)^s$ for some $s=1,2$ and from the
Wronskian we obtain $\vp_{p+1}(\m_n)\vt_{p}(\m_n)=1$, then
$|\vp_{p+1}(\m_n)|=1$.

Let $\m_n\in \g_n$. Due to periodicity we have
$$
\vp_{j+p}(\l)=A\vp_{j}(\l)+B\vt_{j}(\l)
$$
and then at $j=0,1$ we have $B=\vp_{p}(\l)$ and $A=\vp_{p+1}(\l)$.
Thus at $\l=\m_n$ we obtain $\vp_{j+p}(\m_n)=A\vp_{j}(\m_n),
A=\vp_{p+1}(\m_n)$. Thus we have \er{bfc}. \BBox

\subsection{Inverse problem}

We need the following results from the inverse spectral theory
 for the operator $J_+$ on the half-line, in the form
convenient for us. Let $v_{1x}=\log a_x\in R, v_{2x}=b_x$. We can
take $(a_x, b_x)$ as a vector $v$ in the form:
 \[
 \lb{a}
  v=(v_{1x},v_{2x})_1^p\in\mH^2,\ \ \mH:=\lt\{b\in\R^{p}:\ \sum_1^p
b_x=0\rt\}.
\]
Here we have $(v_{1x})_1^p\in \mH$, since $\prod_1^p a_x=1$.
 Using symmetrization, we construct a gap length mapping
 $\p:\mH^2\to \R^{2p-2}$  by
 $$
v\to  \p(v)=(\p_x(v))_1^{p-1},\qqq \p_x=(\p_{1x},\p_{2x})\in\R^2,
$$
and the components have the form
\[
\lb{c}
\begin{aligned}
 \p_{1n}={\l_n^+ + \l_n^-\/2} -\m_n, \qqq
\p_{2n}=\rt|{|\g_n|^2\/4}-\p_{1n}^2\rt|^{1\/2}\e_n, \qq \ca \e_n=1\
\ & {\rm if}\ \m_n^+\in  \g_n^1,
\\
\e_n=-1 \ \ & {\rm if}\ \m_n^+ \in \g_n^{2},
\\
\e_n=0 \ \ & {\rm if}\  \m_n^+ \in \{\l_n^+, \l_n^- \} \ac.
\end{aligned}
\]
Due to \er{bfc} we have  $\e_n=\sign \log |\vp_{p+1}(\m_n)|$ and
note that $\vp_{p+1}(\m_n)\ne 0$.  In order to construct the vector
$\p$ we need:  the gap length $|\g_n|$, $\p_{1n}$  and the sign
$\e_n$ for all $n\in \N_{p-1}$. We formulate the result about the
mapping $\p $, which is similar to results from \cite{K06} and it is
some analogous of the gap-length mapping for periodic Schr\"odiger
operators on $\R_+$ from \cite{K99}.

 \begin{theorem} \lb{Ta}
The mapping $\p:\mH^2\to \R^{2p-2}$ given by \er{c} is a real
analytic isomorphism between $\mH^2$ and $\R^{2p-2}$.

\end{theorem}

{\bf Proof.} In \cite{K06} we consider the mapping,  where $\m_n$
are the zeroes of $\vt_{p+1}(\l)$, i.e. , we use the Neumann
eigenvalues. In the present paper we discuss the case, where  $\m_n$
are the zeroes of $\vp_{p}(\l)$, i.e. , we use the Dirichlet
eigenvalues. We omit the proof of theorem for the Dirichlet
eigenvalues, since it is very similar to the case of the Neumann
eigenvalues in \cite{K06}. \BBox


\section {One dimensional  half-solid}
\setcounter{equation}{0}

We discuss a  half-solid model in $\Z$. In this case we consider the
Jacobi operator  $T_\t, \t\in \R$ on $\ell^2(\Z)$ given by
\[
\lb{b1} (T_\t f)_x=a_{x-1}f_{x-1}+a_{x}f_{x+1}+b_xf_x,\qqq n\in\Z,
\]
where  $\t>1$ is large enough and  the coefficients $a_x,b_x$
satisfy
\[
\lb{ab}
\begin{aligned}
\ca  a_x, b_x, x\ge 1, \ {\rm are \ p-periodic}, \qq \prod_1^p
a_j=1,
\\
a_x=1, b_x=\t,\ x\le 0 \ac .
\end{aligned}
\]
By the physical point of view,  $b_x, x\ge 1$ is a crystal potential
and the real constant $\t$ is the vacuum potential. Define two
operators $J_+$ on $\ell^2(\Z_+)$ and  $J_\t$ on $\ell^2(\Z_-)$ by
\[
\lb{1+2}
\begin{aligned}
 &  (J_+ f)_x= a_{x-1}f_{x-1}+a_{x}f_{x+1}+b_xf_x,\qq x\ge 1, \qq f_0=0,
\\
&  (J_\t f)_x=f_{x-1}+f_{x+1}+\t f_x, \qq x\le 0, \qq f_1=0.
\end{aligned}
\]
Let $P_\pm$ be the projector from $\ell^2(\Z)$ onto
$\ell^2(\Z_\pm)$. We rewrite the operator $T_\t$ in the form
\[
\lb{ab2}
\begin{aligned}
T_\t P_\pm=J_\pm P_\pm,\qqq
   \ca (T_\t f)_1 =a_0f_0+ (J_+ P_+f)_1
\\
 (T_\t f)_{-1} =f_0+ (J_\t  P_-f)_{-1}
 \\
 (T_\t f)_{0} = f_{-1}+ a_0f_{1}+\t f_0
 \ac .
\end{aligned}
\]
In fact, we discuss the case of one-dimensional octant periodic
potentials in the specific form given by \er{ab}. In order to
describe the spectrum of $T_\t$ we use some properties of the
operator $J_+$ on the half-line $\Z_+$ from Section 2. We recall
needed results about  operators $T_\t$.  We have the following
simple results about the spectrum of  $\s (T_\t)$ given by
\[
\lb{Tac}
\begin{aligned}
\s(T_\t)=\s_{ac}(T_\t)\cup \s_{d}(T_\t),\qq
 \s_{ac}(T_\t)=\s_{ac}(J_+)\cup \s_{ac}(J_\t ),\qq \s_{ac}(J_\t )=[\t-2,\t+2].
\end{aligned}
\]
Recall that we assume that the parameter $\t>1$ is large enough. In
this case we have
\[
\begin{aligned}
 \s_{ac} (T_\t)=\cup_{n=0}^p \s_n(T_\t), \ \  \s_n(T_\t)=\s_n(J_+),
n=0,1,..,p-1, \qq \s_p(T_\t)=[\t-2,\t+2].
\end{aligned}
\]
Thus, all possible gaps in the spectrum $\s_{ac}(T_\t)$ are given by
\[
\lb{g1}
\begin{aligned}
\g_n(T_\t)=\g_n(J_+),\qq  n\in \N_{p-1},\qq  \g_p(T_\t)=(\l_p^+,
\t-2).
\end{aligned}
\]

We begin to describe eigenvalues of $T_\t$.  For the operator $T_\t$
we introduce the Jacobi equation
\[
\lb{b2} a_{x-1}f_{x-1}+a_{x}f_{x+1}+b_xf_x=\l f_x,\qqq x\in\Z.
\]
For the operator $J_+$ we define the Weyl  function $\p_x^+$ by
\[
\p_x^+=\vt_x+m_+\vp_x,\qqq  x\ge 1 ,
\]
where $\vt_x, \vp_x$ are solutions of the equation \er{b2} under the
conditions $\vp_{0}=\vt_1=0$ and $\vp_1=\vt_{0}=1$. Note that
$\p_x^+$ depends on $a_p, a_x, b_x, x\ge 1$ only.

For the operator $J_\t $ we define the Weyl  function $\p_x^-, x\le
0$. The equation \er{b2} has the form
\[
\lb{bc} \p_{x-1}^-+\p_{x+1}^-=(\l-\t) \p_x^-=(z+{1\/z})\p_x^-,\qqq
x\le 0,
\]
where $z\in \dD_1=\{\l\in \C: |\l|<1\}$ is defined by
$\l-\t=z+{1\/z}$. Thus we have
\[
\lb{Je}
\begin{aligned}
\p_x^-=z^{-x},  \qq x\le 1,
\\
z=z(\l)=t-\sqrt{t^2-1}\in \dD_1,\qq t={\l-\t\/2}, \\
 z(\l)={1\/2t}+{O(1)\/t^3}\qqq \as \ t\to \iy.
\end{aligned}
\]
For the operator $T_\t$ we introduce the Weyl-type functions
$\P_x^{\pm}(\l)$,  which are solutions of the
 equation \er{b2} and satisfy
$$
(\P_x^{\pm}(\l))_{x\in \Z_\pm}\in \ell^2({\Z_{\pm}}), \qqq \forall \
\l\in \cL:=\C \sm \s_{ac}(T_\t).
$$
For $\l\in \cL$ they  have the forms
\[
\lb{pp}
\begin{aligned}
\P_x^{+}(\l)=\p_x^+(\l)=\vt_x(\l)+m_+\vp_x(\l), \ \  n \ge 1; \\
\P_x^{-}(\l)=z^{-n},  \ \  n \le 1.
\end{aligned}
\]
 These functions $\P_x^{\pm}(\l)$ are analytic in the cut domain
$\cL$ and are continuous up to the boundary. We compute $\P_0^+$.
From \er{b2} and \er{pp}, we get
\[
\lb{11}
\begin{aligned}
\P_{0}^++a_1\p_{2}^++ (b_1-\l)\p_{1}^+=0 \Leftrightarrow \P_{0}^+
-a_p=0 \Leftrightarrow  \P_{0}^+={a_p}.
\end{aligned}
\]
Thus due to \er{pp} -\er{11} and $a_0=1$ we obtain
\[
\lb{wpp}
\begin{aligned}
w=\{\P^-,\P^+\}_0=a_0(\p^-_0\p^+_{1}-\P^+_0\P^-_{1})=m_+-{a_p\/ z}.
\end{aligned}
\]
The function $w(\l)$ is analytic on the domain $\cL$ and has finite
number of zeros, which are simple and coincide with eigenvalues of
the operator $T_\t$. In Lemma \ref{TD1} we show that in each open
gap $\g_j(T_\t)\ne \es, j\in \N_{p-1}$ there is at most one
eigenvalue $\m_j(\t)\sim \m_j$ at large $\t$. We discuss the
eigenvalues of $T_\t$ in the gaps $\g_n(T_\t), n\ge 0,$ and
determine how these eigenvalues depend on $\t$ large enough.

\begin{lemma}
\lb{TD1} Let the operator $J_+$
 on  $\ell^2(\Z_+)$ defined by \er{1+2}   have an
open gap $I=(\l^-,\l^+)$ in the continuous spectrum and  an
eigenvalue $\m\in (\l^-,\l^+)$ for some p-periodic $a,b$. Then for
any constant $\t$ large enough the operator $T_\t$ defined by
\er{b1}, \er{ab} has exactly one eigenvalue $\m_\t$ in the gap $I$
such that
\[
\lb{ms} \m_\t=\m+{c(\m)\/\t}+{O(1)\/\t^2}  \qqq \as \qqq \t\to \iy,
\]
where $c(\m)={2\gF_o(\m)\/a_p\vp_p'(\m)}\ne 0$. Moreover, if $J_+$
has a resonance on the interval $I^2=(\l^-,\l^+)$ on the second
sheet of the operator $J_+$, then for any constant $\t$ large enough
the operator $T_\t$ defined by \er{b1}, \er{ab} has not any
eigenvalue in the gap $I$.
\end{lemma}

\no {\bf Proof.} Using \er{wpp}, \er{Je} we rewrite the Wronskian
$w(\l)$ in the gap $\g_n\ss \cL$ in the form
\[
\lb{w1}
\begin{aligned}
w(\l)=m^+(\l)-{a_p\/z}={\gF_o(\l)-b(\l)\/\vp(1,\l)}-a_pz_1(\l),
\\
z_1(\l)= t+\sqrt{t^2-1}, \qq t={\l-\t\/2}, \qqq \l\in\g_n\ss \cL,
\end{aligned}
\]
since $z z_1=1$, and
\[
\lb{w2}
\begin{aligned}
z_1(\l)=t+\sqrt{t^2-1}={\l-\t}-{O(1)\/\t}, \qq {\as}\  \ \t\to \iy,
\l\in \cL,
\\
 (-1)^nb(\l)=\sqrt{\gF^2(\l
)-1}>0, \qq \qq {\rm if}\  \  \l\in \g_n \ss \cL.
\end{aligned}
\]
  The eigenvalues of $T_\t$ are zeros of  the Wronskian $w(\l)$,
given by \er{w1}, on the domain $\cL$. Consider the two functions
$m_+(\l)$ and $z_1(\l)$ on the gap $(\l^-,\l^+)$, where $\t\to \iy$.
The point $\m\in I$ is an eigenvalue of the operator $J_+$. Then due
to \er{wf} we have $\gF_o^2(\m)=b^2(\m)\ne 0$ and
$\gF_o(\m)=-b(\m)\ne 0$ since the functions $m_+(\l)$ has the pole
at $\m_n\in \g_n^1$. Then the function $m_+(\l)$ is a meromorphic in
the disk $\{\l\in \L_1:|\l-\m|<\ve \}$ around  $\m\in \g_n^1$ and
has the following asymptotics
\[
\begin{aligned}
{m_+(\l)\/a_p}={c(\m)\/\l-\m}+O(1)\qqq \as \qq \l\to \m,\qq
c(\m)={2\gF_o(\m)\/a_p\vp_p'(\m)},
\end{aligned}
\]
We have also $ z_1(\l)= \l-\t+{O(1)\/\t}$ as $\t\to +\iy$ locally
uniformly in $\l\in \C$.
 Thus  the equation
$m_+(\l)=z_1(\l)$ has a unique solution $\m_\t\to \m$ as $\t\to \iy$
given by \er{ms}, since $ {c(\m)\/\m_\t-\m}=\t-\l+O(1). $

Let  $J_+$ have a resonance on the interval $I^2=(\l^-,\l^+)$ on the
second sheet of the operator $J_+$. Then due to Lemma \ref{TAn} the
function $m_+$ is analytic the interval $I=(\l^-,\l^+)$ on the first
sheet of the operator $J_+$ and the function $m_+$ is uniformly
bounded on $[\l^-,\l^+]$. Then due to the simple asymptotics
\er{Je}, the Wronskian $w=m_+-{a_p\/ z}$ has not any zero on
$I=(\l^-,\l^+)$ for any constant $\t$ large enough. Thus  for any
constant $\t$ large enough the operator $T_\t$ defined by \er{b1},
\er{ab} has not any eigenvalue in the gap $I$. \BBox

Now we prove the main result of this section. Recall that
$\N_m=\{1,2,...,m\}$.

\begin{lemma}
\lb{TD2} i) Let integer $p\ge 2$ and let $\g>0$. Then there exist
p-periodic sequences $a_n, b_n$ such that all $p-1$ gaps in the
spectrum of the operator $J_+$ on $\ell^2(\Z_+)$ are open and
satisfy
\[
|\g_j|=\g,\qq \forall \ j\in \N_{p-1}.
\]
In addition, for any points $\l_j\in\g_j^c, j\in\N_{p-1} $, exist
unique p-periodic sequences $a_n, b_n$ such that each $\l_j=\m_j^+,
$ is a state of the operator $J_+$.

ii) Let in addition the operator $T_\t$ be given by \er{b1}
 \er{ab} and let $\t$ be large enough. If $\m_j^+\in \g_j^1$ is an
 eigenvalue of the operator $J_+$, then the operator  $T_\t$
 has a unique eigenvalue $\m_j(\t)$ on the gap $\g_j^1$ such that
 for some constant $c(\m)\ne 0$:
\[
\lb{mns} \m_j(\t)-\m_j^+={c(\m_j)\/\t}+{O(1)\/\t^2}  \qqq \as \qqq
\t\to \iy.
\]

If $\m_j^+\in \g_j^2$ is a resonance of the operator $J_+$ for some
$\in \N_{p-1}$, then the operator  $T_\t$ has not eigenvalues on the
gap $\g_j^1$.

\end{lemma}

\no {\bf Proof} of i)  follows from  Theorem \ref{Ta}. The proof of
ii) follows from i) and Lemma \ref{TD1}.
 \BBox

\

\section {\lb{Sec3} Difference operators on the lattice}
\setcounter{equation}{0}

\subsection{Specific periodic  Jacobi  operators on the half-line}
Consider the  Jacobi   operator $J_+$ on $\ell^2(\Z_+)$ given by
\er{Jdef}.  Recall that the spectrum of $J_+$ consists of an
absolutely continuous part  (which is a union of non-degenerate
spectral bands $\s_n=[\l^+_n,\l^-_{n+1}], n=0,1,..,p-1$) plus at
most one eigenvalue in each open  gap $\g_n=(\l^-_{n},\l^+_n), n\in
\N_{p-1}$.

{\it Now we begin to construct a specific Jacobi   operator $J_+$.}
 Here we
use results about the gap-lengths mapping from Lemma \ref{TD2} i).
Due to these results about the gap-lengths mapping,
 we take the coefficients  $a_n,b_n$ such that all $p-1$ gaps
 $\g_1, ..., \g_{p-1}$ are open in the spectrum of $J_+$  and  satisfy
\[
\lb{bg1}
\begin{aligned}
\l_0^+=0,\qqq
 \g=|\g_1|=|\g_2|=|\g_3|=....|\g_{p-1}|.
\end{aligned}
\]
Let $\gS_n=\sum_{j=0}^{n-1}|\s_j|$. Thus \er{bg1} and the estimate
\er{eb} give
\[
\lb{bg1A}
\begin{aligned}
\l_n^-=\g(n-1)+\gS_n, \qqq \l_n^+=\g n+\gS_n, \qqq  |\gS_{n}|\le 4.
\end{aligned}
\]
 Due to Lemma \ref{TD2} i) in each gap $\g_n$, $n\in \N_{p-1}$ of $J_+$ we choose exactly one
eigenvalue $\m_n^+$ by
\[
\lb{bg2} \m_n^+=\g e_n \in \g_n^1,\qqq e_n=n-1+e_1,\qq e_1={1\/4d}.
\]
It is convenient to define the {\bf normalized} operator $
\cJ_\g={1\/\g}J_+$.  Then the spectrum of $\cJ_\g$ consists of union
of bands $s_n={1\/\g}\s_n$  part plus exactly one eigenvalue
$e_n={\m_n^+\/\g}$ in each open gap $g_n={1\/\g}\g_n$. Thus due to
\er{bg1}, \er{bg1A} we have
\[
\lb{ebx}
\begin{aligned}
\textstyle & s_0={\s_0\/\g},\qqq
s_n={\s_n\/\g}=[n+{\gS_n\/\g},n+{\gS_{n+1}\/\g}], \qq
g_{n}={\g_n\/\g},\qq
\\
& |s_0|\le {4\/\g},\qqq  |s_n|\le {4\/\g}, \qqq |g_n|=1,
\end{aligned}
\]
for all $n\in \N_{p-1}$,  where $\gS_n$ is defined in \er{bg1A}.
Thus each  spectral band $s_n$ is very small and is very close to
the point $n$ and satisfies
\[
\lb{snn} \dist \{s_n,n\}\le {\gS_{n+1}\/\g}\le {4\/\g}, \qq
n=0,1,...,p-1.
\]
In each open gap $g_n, n\in \N_{p-1}$, there exists exactly one
eigenvalue $e_n$ of $J_\g$ such that
\[
\lb{es2}
\begin{aligned}
\textstyle
 e_n={\m_n\/\g}.
\end{aligned}
\]


\subsection{ Difference operators on $\Z_+^2$}  We consider  the difference periodic
operator $H_0=J_1+J_2$ on the corner $\Z_+^2$ acting on the
functions $f_x, x=(x_1,x_2)\in \Z_+^2$. Here $J_j,j=1,2$ is the  p
periodic Jacobi operator on the half-lattice $\Z_+$ and given by
$$
(J_j f)_{m}=a_{m-1}f_{m-1}+a_{m}f_{m+1}+b_{m}f_{m},\qqq f_0=0,\qq
m=x_j\in \Z_+=\{1,2,...\}.
$$
We assume that the Jacobi operators $J_1$ and $J_2$ satisfy
\er{bg1}-\er{bg2} with large  gaps in the spectrum.
 For a large constant $\g$ we define a new {\bf normalized} operator
$ \cH_\g={H_0\/\g}={J_1+J_2\/\g}. $ We take the operator $\cH_\g$,
when the variables are separated.  We show that $\cH_\g$ has bands
which are very small and their positions are very close to the
integer $n$. The union of group of bands close to the integer $n$
forms a cluster. Between the two neighbor clusters there exists a
big gap. On this gap there exist eigenvalues. We describe these
clusters and eigenvalues.

$\bu$ We define the basic bands $S_{i,j}^0$ of the operator
${\cH_\g}$ and their clusters $K_n^0$ by
\[
\lb{Ks} S_{i,j}^0=s_i+s_j,\qqq i,j=0,1,....,p-1,\qqq
K_n^0=\bigcup_{i+j=n} S_{i,j}^0,\qqq n=0,1,....,2p-2,
\]
where we define $A+B$ for sets $A,B$ by $ A+B=\{z=x+y: (x,y)\in A\ts
B\}$. In particular, we have
\[
\begin{aligned}
\lb{K1} K_0^0=S_{0,0}^0,\qq K_1^0=S_{1,0}^0,\qq K_2^0=S_{2,0}^0\cup
S_{1,1}^0,.....,
\end{aligned}
\]
If $\g$ is large enough, then due to \er{ebx}, \er{snn}  we estimate
the position of bands $S_{i,j}^0$,  their lengths $|S_{i,j}^0|$ and
their cluster $K_n^0$  by
\[
\lb{K2} \dist \{S_{i,j}^0,i+j\}\le {2\/\g},\qqq \dist \{K_n^0,n\}\le
{2\/\g}.
\]
$\bu$ A surface band is created by an eigenvalue $e_j$ and  a band
$s_i$ of Jacobi operators.  We define the surface bands $S_{i,j}^1$
and their clusters $K_n^1$ of the operator ${\cH_\g}$ by
\[
\begin{aligned}
S_{i,j}^1=e_i+s_j,\qq i\in \N_{p-1}, \qq j=0,1,2,3,...,
%
\qq
 K_n^1=\bigcup_{i+j=n+1} S_{i,j}^1,\qqq , n=0,1,....,p.
\end{aligned}
\]
In particular, we have
\[
\begin{aligned}
 K_0^1=S_{1,0}^1,\qq
K_1=S_{2,0}^1\cup S_{1,1}^1,\qq
 K_2^1=S_{3,0}^1\cup S_{1,2}^1\cup S_{2,1}^1,\qq .....
\end{aligned}
\]
The position of the guided bands $S_{i,j}^1$ and the cluster $K_n^1$
is given by
\[
S_{i,j}^1\sim (i+j-1)+e_1=n+e_1,\qqq K_n^1\sim n+e_1.
\]

\no  $\bu$ The operator $\cH_\g$ has {\bf eigenvalues} $K_n^e, n\ge
0 $ with multiplicity n+1  given by
\[
\lb{Ke1} K_n^e=E_{i,j}:=e_i+e_j=n+2e_1,\qqq e_1={1\/4d},\qqq \
i+j=n+2,n \ge 0
\]
for all $i,j=1,2,...,p$. In particular, we have
\[
\begin{aligned}
\lb{Ke1x} K_0^e=E_{1,1},\qqq K_1^e=E_{1,2}=E_{2,1},\qqq
K_3^e=E_{1,3}=E_{2,2}=E_{3,1},....
\end{aligned}
\]
$\bu$ Thus we can describe  $\s_{ac}(\cH_\g)$ and
$\s_{disc}(\cH_\g)$ by
\[
\s_{ac}(\cH_\g)=\cup_{n\ge 0} (K_n^0\cup K_n^1),\qqq
\s_{disc}(\cH_\g)=\cup_{n\ge 0}K_n^e,
\]
where
\[
\lb{k12d} K_n^0\sim n,\qq K_n^1\sim n+e_1, \qqq K_n^e=n+2e_1.
\]
 We have two types of band clusters  $K_j^0$ and $K_j^1$. These clusters  are separated by
gaps. Now  combining all estimates \er{Ks}-\er{k12d} we deduce that
there exists an interval $I_n$ such that
\[
\lb{Inx}
\begin{aligned}
I_n=[K_n^e-r,K_n^e+r],\qqq \s_{ac}(\cH_\g)\cap I_n=\es, \qqq {\rm
where}\ r={e_1\/2}={1\/8d},
\end{aligned}
\]
  for some $\g>0$ large enough. Thus the
spectral interval $\gI_{n,\g}=\g I_n$ satisfies
\[
\lb{Inxx}
\begin{aligned}
   \gI_{n,\g}=\g I_n=[\g (K_n^e-r), \g (K_n^e+r)],
   \qqq
   \dist\{\gI_{n,\g}, \s_{ac}(\cH_0)  \}\ge 3.
\end{aligned}
\]
Then  interval $\gI_{n,\g}\cap \s_{ac}(H_0)=\es$ and the operator
$H_0$ has the eigenvalue $\g K_n^e\in \gI_{n,\g}$ of multiplicity
$n+1$. Moreover, the interval $\gI_{n,\g}$ does not contain other
spectrum and  to the right and to the left of it there is a
essential spectrum. In fact we have proved Theorem \ref{T1} ii) for
the case $H_0$.

\subsection{\lb{H3} Difference operators on $\Z_+^3$}
We consider  difference operators $J=J_1+J_2+J_3$ on the corner
$\Z_+^3$  and acting on the functions $f_x, x=(x_1,x_2,x_3)\in
\Z_+^3$. Here $J_j,j=1,2,3$ is the  p periodic Jacobi operator on
the half-line $\Z_+$ and given by
$$
(J_j f)_{m}=a_{m-1}f_{m-1}+a_{m}f_{m+1}+b_{m}f_{m},\qqq f_0=0,\qq
m=x_j\in \Z_+.
$$
We assume that the Jacobi operators $J_j$ satisfy \er{bg1}-\er{bg2}
with large  gaps in the spectrum.
 For large constant $\g$ we define a new {\bf normalized} operator
 by
$$
\cJ_\g={\cJ\/\g}={J_1+J_2+J_3\/\g}=J_{1,\g}+J_{2,\g}+J_{3,\g},\qqq
J_{j,\g}={J_j\/\g}.
$$
$\bu$ We define basic bands $S_{i,j,k}^0$ of the operator ${\cJ_\g}$
and their clusters $K_n^0, n=0,1,....$ by
\[
S_{i,j,k}^0=s_i+s_j+s_k,\qqq i,j=0,1,2,.., \qq k\in \N_{p-1},\qqq
K_n^0=\bigcup_{i+j+k=n} S_{i,j,k}^0,
\]
and in particular,
$$
\begin{aligned}
& K_0^0=S_{0,0,0}^0=s_0+s_0+s_0,\qqq K_1^0=S_{0,0,1}^0,\qqq
K_2^0=S_{0,0,2}^0\cup S_{0,1,1}^0,....
\end{aligned}
$$
Recall that  we define $A+B$ for sets $A,B$ by $A+B=\{z=x+y:
(x,y)\in A\ts B\} $. Similar to 2-dim case  we deduce that
\[
S_{i,j,k}^0\sim i+j+k, \qqq K_n^0 \sim n,\qqq  \forall\qq
n=1,2,...,N.
\]

In 3-dimensional case we have two types of the surface bands
$S_{i,j,k}^1$ and $S_{i,j,k}^2$.

\no $\bu$  {\bf The first  type of surface bands.} We define the
surface bands $S_{i,j,k}^1$ of the operator ${\cJ_\g}$ and their
clusters $K_n^1, n=1,2,....$ by
\[
S_{i,j,k}^1=s_i+s_j+e_k,\qqq  K_n^1=\bigcup_{i+j+k=n+1}
S_{i,j,k}^1,\qqq i,j=0,1,2,.., \qq k\in \N_{p-1}.
\]
The position of surface bands $S_{i,j,k}^1$ and their clusters
$K_n^1$ are given by
\[
S_{i,j,k}^1\sim i+j+k-1+e_1=n-1+e_1, \qqq \qqq K_n^1\sim n-1+e_1.
\]
These clusters  are separated by gaps.  Thus we have
\[
\begin{aligned}
K_1^1=S_{0,0,1}^1,\qq K_2^1=S_{0,0,2}^1\cup S_{0,1,1}^1,\qq
K_3^1=S_{0,0,3}^1\cup S_{0,1,2}^1\cup S_{1,1,1}^1,....
\end{aligned}
\]
$\bu$ {\bf The second type of surface (guided) bands.} We define the
surface (guided) bands $S_{i,j,k}^2$ of the operator ${\cJ_\g}$ and
their clusters $K_n^2, n=0,1,....$ by
\[
S_{i,j,k}^2=e_i+e_j+s_k,\qqq \qqq K_n^2=\bigcup_{i+j+k=n+2}
S_{i,j,k}^2,\qq i,j\in \N_{p-1},\qq k=0,1,2,...
\]
The positions of the surface bands $S_{i,j,k}^2$ and the cluster
$K_n^2$ are given by
\[
S_{i,j,k}^2\sim i+j-2+2e_1+k=n+2e_1, \qqq K_n^2\sim n+2e_1.
\]
These clusters  are separated by gaps.  Thus we have
\[
\begin{aligned}
K_0^2=S_{1,1,0}^2,\qq K_1^2=S_{1,1,1}^2\cup S_{2,1,0}^2,\qq
K_2^2=S_{3,1,0}^2\cup S_{2,1,1}^2\cup S_{1,1,2}^2,.....
\end{aligned}
\]
$\bu$ {\bf Eigenvalues.} Due to \er{es2} the operator $\cJ_\g$ has
eigenvalues $K_n^e$ given by
\[
K_n^e=e_i+e_j+e_k=i+j+k-3+3e_1=n+3e_1,\qq i,j,k\in \N_{p-1},\qq
i+j+k=n+3,
\]
$n=0,1,..$. The sets $\s_{ac}({\cJ_\g})$ and $\s_{disc}({\cJ_\g})$
are given by
\[
 \s_{ac}(\cJ_\g)=\cup_{n\ge 0} (K_n^0\cup K_n^1\cup K_n^2),\qqq
\s_{disc}(\cJ_\g)=\cup_{n\ge 1}K_n^e.
\]
Later on we repeat the consideration for the case $d=2$.

\subsection{\lb{h} Specific 1dim half-solid  potentials}
 Consider the  Jacobi  operator as a half-solid model in $\Z$. In this case we
consider the Jacobi operator  $T_\t$ on $\ell^2(\Z)$ given by
\[
\lb{b1x} (T_\t f)_x=a_{x-1}f_{x-1}+a_{x}f_{x+1}+b_xf_x,\qqq x\in\Z.
\]
Let  $\t$ be large enough and  the coefficients $a_x,b_x$ satisfy
\[
\lb{abx}
\begin{aligned}
\ca  a_0, a_x, b_x,  \ {\rm are \ p-periodic}\qq & x\ge 1,
\\
a_x=1, b_x=\t=b_0, & x\le -1 \ac
\end{aligned}
\]
Let an integer $p\ge 2$ be large enough. Due to Lemma \ref{TD2} for
large $\g>1$  we obtain that there exists p-periodic $a_x, b_x, x\ge
1$ sequences such that \er{bg1} holds true. Thus  by
\er{Tac}-\er{g1}, all gaps $\g_j, j\in \N_{p-1}$ in the spectrum of
the operators $J_+$ and $T_\t$ are open. Moreover, there exists an
eigenvalue $\wt\m_j(\t)$ of $T_\t$ in each this gap $\g_j$ and they
satisfy
\[
\lb{gjxx} \wt\m_j(\t)\in \g_j \qqq |\g_j|=\g,\qq \forall \ j\in
\N_{p-1},
\]
\[
\lb{Ts1} \s_{ac}(T_\t)=\s_{ac}(J_+)\cup \wt\s,\qqq \s_{ac}(J_+)=
\s_0\cup\s_1\cup.... \s_{p-1}, \qqq \wt\s= [\t-2,\t+2].
\]
Here the bands $\s_0, \s_1,.... \s_{p-1}$ are separated by gaps
$\g_j,  j\in \N_{p-1}$ and the bands $\s_{p-1}$ and $\wt\s$ are
separated by a gap $\wt\g_p=(\l_p^+,\t-2) $  and each eigenvalue
$\m_n(\t)$ satisfies \er{mns}.

 Define  a new {\bf normalized} operator
$ T_{\t,\g}={1\/\g}T_\t $. From the properties of $T_\t$ we deduce
that the spectrum of $T_{\t,\g}$ consists of an absolutely
continuous part
\[
\lb{Ts1x} \s_{ac}(T_{\t,\g})= \bigcup\limits_{n=0}^{p} s_n,\qqq
s_n={\s_n\/\g},\qq n\in \N_{p-1},\qq s_p={\wt\s\/\g},
\]
 plus  at
most one eigenvalue in each non-empty finite gap $g_n$, $n\in\N_p$,
 given by
$$
 g_{n}={\g_n\/\g}, \qq n\in \N_{p-1},
$$
and they satisfy \er{ebx}-\er{es2}. In each gap $g_n, n\in
\N_{p-1}$, there exists exactly one eigenvalue $\wt e_n$ given by
\[
\lb{enx} \wt e_n={\wt\m_n\/\g}=e_n+\ve_n,\qqq e_n=n-1+{1\/4d},\qq
|\ve_n|\le {1\/\g},\qqq    n\in \N_{p-1},
\]
since $\t>1$ is large enough. Thus roughly speaking the spectrum of
the operators $T_{\t,\g}$ on $\ell^2(\Z)$ and $\cJ_\g$ ( on
$\ell^2(\Z_+)$) is the same on the interval $[0,\l_p^+]$. They have
the same bands $\s_0, \s_j, j\in \N_{p-1}$ and the same gaps $\g_j,
j\in \N_{p-1}$. Moreover, their eigenvalues $\wt e_j$ and $e_j$ in
each gap $\g_j$ are very close, since we take $\t$ large enough.


\subsection{\lb{H2R} Model difference operators on $\Z^2$}
We consider difference operators $H_0=T_{\t,1}+T_{\t,2}$ on the
lattice  $\Z^2$,   where $T_{\t,1},j=1,2$ is the Jacobi operator on
the lattice $\Z$, discussed in Subsection 4.4.  The spectrum of
$T_{\t,j}$ and $J_{\t,j}$ are similar on the interval
$[\l_0^+,\l_{p}^+]$. Then the spectrum of the sum
$T_{\t,1}+T_{\t,2}$ is similar to the spectrum of $J_{1}+J_{2}$ on
the interval $[\l_0^+,2\l_p^+]$. The proof repeats the case
$J_{1}+J_{2}$. Moreover, using similar arguments
 we prove Theorem \ref{T1} for the operator $T=T_{\t,1}+T_{\t,2}$.
 The proof for the case $\Z^d, d\ge 3$ is similar.

\subsection{Model difference operators on $\Z_+\ts \Z$}
Consider the operator $H_0=J_{1}+T_{\t,2}$ on the half-lattice
$\Z_+\ts \Z$, where the operator $J_{1}$ acts on the half-line and
depends on one variable $x_1\in \Z_+$; the operator $T_{\t,2}$
(depending on one variable $x_2\in \Z$) acts on $\Z$ and given by
\er{b1x}, \er{abx} and the constant $\t\ge 1$ is large enough. The
spectrum of $T_{\t,2}$ and $J_1$ are similar on the interval
$[\l_0^+,\l_p^+]$ for $p, \t$ large enough. The proof repeats the
case $J_1+J_2$. Moreover, using similar arguments
 we prove Theorem \ref{T1} for the operator $H_0=J_{1}+T_{\t,2}$.

\section {Proof of main Theorems}
\setcounter{equation}{0}

{\bf Proof Theorem \ref{T1}} i) We consider an operator $H_+=-\D+V$
on $\R_+^2$, where the potential $V\in \ell^\iy(\Z_+^2)$ is octant
periodic, the proof for other cases is similar. Without loss of
generality, we assume that $V$ is $(p\Z)^2$ -periodic for some
$p>1$. Let $H=-\D+V$ on $\Z^2$. Define functions $g_n\in
\ell^\iy(\Z_+)$ and $G_n\in \ell^\iy(\Z_+^2), n\ge 1$ by:
\[
\lb{gn}
\begin{aligned}
g_n|_{w_n}=1, \qqq g_n|_{\Z_+\sm w_n}=0,\qqq w_n=[4^n,4^n+n+1],\qqq w_n\cap w_{n+1}=\es,\\
G_n(x)=g_n(x_1)g_n(x_2),\qq x=(x_1,x_2)\in \Z^2,\qqq \supp G_n\ss
\Z_+^2.
\end{aligned}
\]
Let $\cT=\Z^2/(p\Z)^2$. For any $\l\in \s(H)$ there exists a
function $\p_x=e^{i(k,x)}u(x,k)$, which  satisfies
\[
\lb{pnu}
\begin{aligned}
(-\D +V(x))\p_x(k)=\l\p_x(k),\qq \forall \ x\in \Z^2,\\
u(\cdot,k)\in \ell^2(\cT),\qqq \sum_{x\in\cT}|u(x,k)|^2=1,
\end{aligned}
\]
see \cite{GKT93} for some $k\in \R^2$. For this fix  $k\in \R^2$ we
define the sequence $f_n(x)={1\/c_n}G_n(x)\p_x$, where $c_n>0$ is
given by
$$
c_n^2=\sum_{x\in\Z_+^2}|G_n(x)\p_x(k)|^2.
$$

The function $u(x,k)$ is $(p\Z)^2$ periodic, then due to \er{pnu} we
obtain
\[
c_n^2=\sum_{x\in\Z^2}|G_n(x)\p_x(k)|^2=n^2 +O(n)
\]
as $n\to \iy$. Thus the sequence $f_n$ satisfies

1) $\|f_n(\cdot,k)\|=1$ and $\D f_n\in \ell^2(\Z_+^2)$, for all
$n\in \N$,

2) $f_n \perp f_m$ for  all $n\ne m$, and $f_n\to 0$ weakly as $n\to
\iy$.

Thus $\l\in \s_{ess}(H_+)$, which yields \er{spH}, since standard
arguments imply
$$
\|(H_+-\l)f_n\|=\|(H-\l)f_n\|\to 0\qq  as \qq n\to \iy.
$$


We prove ii) for the case $D=\Z_+^2$ and $N\ge 1$, the proof of
other cases is similar. Consider the operator $H_0=J_1+J_2$, where
$J_1+J_2$ is the described in Subsection 4.2 and $J_1, J_2$ are the
Jacobi operator on $\Z_+$. We assume that they have the properties
in \er{bg1}-\er{snn} for some p-periodic sequences $(a_n, b_n)\in
\R_+\ts \R $. Due to \er{Inxx} for each $n$ the operator $H_0$ has
the eigenvalue $E=\g(n+2e_1)$ of multiplicity $n+1$ and the the
interval $\gI_{n,\g}$ such that
$$
\gI_{n,\g}=\g I_n=[E-\g r, E+\g r],\qqq  \gI_{\g,r}\cap
\s_{ac}(H_0)=\es, \qqq {\rm where}\qq r={e_1\/2}.
$$
  Moreover, the interval $\gI_{n,\g}$ does not contain
other spectrum  and  to the right and to the left of it there is a
essential spectrum. In fact we have proved ii) for the case $H_0$.

We consider an operator  $H_\ve =H_0+\ve W$ on $\ell^2(\R_+^2)$ and
$W$ is the multiplication operator. Here $H_\ve $ is the difference
operator on the quadrant  $ \Z_+^{2}$ with the Dirichlet boundary
conditions on the boundary $\pa \Z_+^{2}$ with octant periodic
coefficients. We assume that the perturbation $W$ satisfies
\[
\lb{s3x}
\begin{aligned}
W=   \sum_{i=1}^{2} (\wt a^i U_i + U_{-i}\wt a^i)+\wt V,
\\
\|\wt V\|_{\ell^\iy(\Z_+^2)}+\|\wt a^1\|_{\ell^\iy(\Z_+^2)} +\|\wt
a^2\|_{\ell^\iy(\Z_+^2)}\le 1,
\end{aligned}
\]
where $(U_if)_x=f_{x+e_i}$ and $(U_i^*f)_x=f_{x-e_i}$ for
$f=(f_x)\in \ell^2(\Z_+^{2})$ and $i=1,2$. We
also assume that $\wt a^i>0$ and $\wt V$ are the octant periodic
functions on $\Z_+^2$. Thus we obtain
\[
\lb{s3} \qq \|W\|\le 5.
\]

We define  contours $c=\{\l\in \C: |\l-E|=2 \}$. Due to \er{Inx} the
operator $H_0$ has     an eigenvalue $E=n+2e_1\in \gI_{\g,r}$ of
multiplicity $n+1$ inside the contours $c$.  Using the simple
identities we deduce that the resolvents $R_0(\z)=(H_0-\z)^{-1}$ and
$R_\ve(\z)=(H_\ve-\z)^{-1}$  satisfy
\[
\begin{aligned}
\| R_\ve(\z)\| \le 2,\qqq \| R_0(\z)\| \le 1,
 \\
\|R_\ve(\z)-R_0(\z)\|\le 5\ve \|R_0(\z)\| \|R_\ve(\z)\|\le 10\ve
\end{aligned}
\]
for all $ \z\in c$ since $\dist \{\s(H_\ve), E\pm 1\}\ge 1-5\ve\ge
{1\/2}$. Then we obtain
\[
\begin{aligned}
P(\ve)=-{1\/2\pi i}\int_{c}R_\ve(\z)d\z,\qq
 \\
R_\ve(\z)=R_0(\z)-R_0(\z)\ve W R_\ve(\z),
\end{aligned}
\]
which yields $ \|P(\ve)-P_n(0)\|<1$. Thus the projectors  $P(\ve)$
and $P(0)$ have the same dimension $n+1$ for all $\ve >0$ small
enough.
 We use similar arguments  in order to show that to the right
and to the left of the interval $[E-1, E+1]$  there is spectra and
its corresponding subspaces have infinite dimension. \BBox

\medskip

\footnotesize \no {\bf Acknowledgments.} \footnotesize This work was
supported by the RSF grant  No. 18-11-00032. We thank Natalia
Saburova for Fig. 2.

\end{document}